\documentclass{amsart}

\usepackage{amsmath}
\usepackage{amsfonts}
\usepackage{amssymb}
\usepackage{amsthm}
\usepackage{graphicx}
\newcommand{\norm}[1]{\left\lVert {#1} \right\rVert}

\newcommand{\C}{{\mathbb{C}}}

\newcommand{\Z}{{\mathbb{Z}}}
\newcommand{\N}{{\mathbb{N}}}

\newcommand{\Ass}{\mathop{\rm Ass}\nolimits}
\newcommand{\ann}{\mathop{\rm ann}\nolimits}

\newtheorem{theorem}{Theorem}
\newtheorem{conj}{Conjecture}
\newtheorem{prop}{Proposition}

\newtheorem{lemma}{Lemma}
\newtheorem{question}{Question}

\theoremstyle{definition}
\newtheorem{defn}{Definition}
\newtheorem{example}{Example}

\theoremstyle{remark}

\author{Jennifer Brooks}
\address{Brooks: Department of Mathematics, Brigham Young University,
Provo,UT, 84602 USA}
\email{jbrooks@mathematics.byu.edu}
\date{\today}

\author{Dusty Grundmeier}
\address{Grundmeier: Mathematics Department, Harvard University, Cambridge, MA, 02138 USA}
\email{deg@math.harvard.edu}

\author{Hal Schenck}
\thanks{Schenck supported by NSF 2006410}
\address{Schenck: Mathematics Department, Auburn University}
\email{hks0015@auburn.edu}


\begin{document}

\title{Algebraic properties of Hermitian sums of squares, II}

\vspace*{-1cm}

\maketitle

\begin{abstract}

We study real bihomogeneous polynomials $r(z,\bar{z})$ in $n$ complex variables for which $r(z,\bar{z}) \norm{z}^2$ is the squared norm of a holomorphic polynomial mapping. Such polynomials are the focus of the Sum of Squares Conjecture, which describes the possible ranks for the squared norm $r(z,\bar{z}) \norm{z}^2$ and has important implications for the study of proper holomorphic mappings between balls in complex Euclidean spaces of different dimension. Questions about the possible signatures for $r(z,\bar{z})$ and the rank of $r(z,\bar{z}) \norm{z}^2$ can be reformulated as questions about polynomial ideals.  We take this approach and apply purely algebraic tools to obtain constraints on the signature of $r$.
\end{abstract}

\section{Introduction}
%
%

Let $r(z,\bar{z})=\sum_{|\alpha| = |\beta| = m} c_{\alpha,\beta}z^\alpha \bar{z}^\beta$ be a real bihomogeneous polynomial on the diagonal of  $\mathbb{C}^n \times \mathbb{C}^n$ of bi-degree $(m,m)$.   Its {\it rank} and {\it signature} are the rank and signature of the corresponding Hermitian coefficient matrix $(c_{\alpha,\beta})$. Suppose the bihomogeneous polynomial $r(z,\bar{z})\norm{z}^2=r(z,\bar{z}) \sum_{j=1}^n |z_j|^2$ is a squared norm, so that
\vskip -.03in
\[
r(z,\bar{z}) \norm{z}^2 = \norm{h(z)}^2=\sum_{k=1}^\rho |h_k(z)|^2
\]
for some holomorphic polynomial mapping $h=(h_1,\ldots,h_\rho)$ defined on $\mathbb{C}^n$.
In this case, the coefficient matrix for $r(z,\bar{z}) \norm{z}^2$ is non-negative semi-definite.

Due to its connection to the Gap Conjecture Problem for proper holomorphic mappings between balls in complex Euclidean spaces of different dimensions, an important open problem in the theory of functions of several complex variables is to determine the possible ranks of $r(z,\bar{z})\norm{z}^2$. 
\begin{conj}[Ebenfelt Sum of Squares (SOS) Conjecture, \cite{E:17}] \label{thm:sos}
Suppose $n \geq 2$.
Let $r(z,\bar{z})$ be a real bihomogeneous polynomial, and suppose
\begin{equation*}
r(z,\bar{z})\norm{z}^2=\norm{h(z)}^2
\end{equation*}
for some holomorphic polynomial mapping $h$. Let $\rho$ be the rank of $\norm{h}^2$  and let
\begin{equation*}
k_0 = \max\left\{k \in \N: \frac{k(k+1)}{2} < n-1\right\}.
\end{equation*}
Then either
\begin{equation}\label{eq:lower bound for rho if N not 0}
\rho \geq (k_0 +1)n -\frac{k_0(k_0+1)}{2},
\end{equation}
or there exists an integer $0\leq k \leq k_0 <n$ such that
\begin{equation}\label{eq:gaps in rank}
n k -\frac{k(k-1)}{2} \leq \rho \leq n k.
\end{equation}
\end{conj}
In \cite{GB}, the first two authors prove the Sum of Squares Conjecture when $r$ has signature $(P,0)$, and \cite{BG:20} shows it holds when $n=3$ and the coefficient matrix of $r$ is diagonal.  The gaps in possible rank described in \eqref{eq:gaps in rank} occur in the non-negative semi-definite case, and it appears that when the coefficient matrix of $r$ has negative eigenvalues, the inequality \eqref{eq:lower bound for rho if N not 0} always holds.

Thus the task is to determine the possible ranks for $r(z,\bar{z}) \norm{z}^2$ when $r$ has signature $(P,N)$ with $N > 0$ and its coefficient matrix is not diagonal.  
Even for $N=1$, we encounter difficulties, leading to the following:
\begin{question}\label{ques:min P for signature P,1}
Suppose $r(z,\bar{z})$ has signature $(P,1)$ and $r(z,\bar{z}) \norm{z}^2$ is a squared norm. What is the minimum possible value for $P$?
\end{question}
When the coefficient matrix for $r$ is diagonal,  \cite{HL:sig} proves
that $P \geq n$ and that $P=n$ is possible.  So a natural question is
whether the same bound holds for arbitrary $r$. Our main result answers this question:
\begin{theorem}\label{thm:sig (2,1) impossible - analysis version}
Let $r(z,\bar{z})$ be a real bihomogeneous polynomial on the diagonal of $\C^n \times \C^n$ of bi-degree $(m,m)$ with signature $(P,1)$.  Suppose that $r(z,\bar{z}) \norm{z}^2$ is a squared norm.  If $n \le 3$, then $P \geq n$. However, if $n \geq 4$, then $P<n$ is possible.
\end{theorem}

\subsection{Algebraic Formulation}

In \cite{GB} and \cite{BG:20}, the first two authors investigate $r(z,\bar{z})\norm{z}^2$ by translating the problem into one about homogeneous ideals in $R=\C[z_1,\ldots,z_n]$. The process begins with a {\it holomorphic decomposition} of $r$; we write  $r(z,\bar{z})=\norm{f(z)}^2 - \norm{g(z)}^2$ where $f \colon \mathbb{C}^n \to \mathbb{C}^P$ and $g \colon \mathbb{C}^n \to \mathbb{C}^N$ are holomorphic polynomial mappings, all of whose components are homogeneous polynomials on $\mathbb{C}^n$ of degree $m$.  We refer to \cite{JPD:Carus} for further discussion of this technique and its applications.

A holomorphic decomposition of a bihomogeneous polynomial is not unique, but if we require the polynomials $f_j,g_k$ to be linearly independent, then $P$ and $N$ are uniquely determined, $(P,N)$ is the signature of the coefficient matrix of $r$, and $P+N$ is its rank.

For $r(z,\bar{z})  = \norm{f(z)}^2 - \norm{g(z)}^2$, consider two ideals:
$I^+$ is the ideal generated by the $P$ components of $f$ and $I^-$ is the ideal generated by the $N$ components of $g$. The condition that $r(z,\bar{z}) \norm{z}^2$ is a squared norm can be restated as a condition on $I^+_{m+1}$ and $I^-_{m+1}$ (the components of these ideals in degree $m+1$), as described in the next proposition. This result first appeared in a somewhat more general form in \cite{JPD:05}; \cite{BG:19} reformulates it in algebraic terms as follows:
\begin{prop}\label{prop:ideal condition for squared norm}
Let $r(z,\bar{z}) = \norm{f(z)}^2 - \norm{g(z)}^2$ be a real bihomogeneous polynomial of bi-degree $(m,m)$, and define $I^+$ and $I^-$ as above.
If $r(z,\bar{z})\norm{z}^2$ is a squared norm,  then 
 $I^-_{m+1} \subseteq I^+_{m+1}$.
\end{prop}
\noindent In this language, Theorem~\ref{thm:sig (2,1) impossible - analysis version} takes the form:
\begin{theorem}\label{thm:sig (2,1) impossible - alg version}
Let $I^+=\langle f_1,\ldots,f_P \rangle$ and $I^- = \langle g \rangle$ for homogeneous polynomials  $f_j,g$ in $R$ of degree $m$, and suppose $\{f_1,\ldots,f_P,g\}$ is a linearly independent set. If $I^-_{m+1} \subseteq I^+_{m+1}$, then $n\le 3 \mbox{ implies } P \geq n$, while if $n\geq 4$, then $P<n$ is possible.
\end{theorem}

\noindent {\bf Acknowledgements.} The first two authors would like to thank several colleagues for patiently entertaining the algebra questions of a pair of complex analysts during the early phases of this work.   Thanks to Steve Humphries, Jason McCullough, and Artan Sheshmani.

\pagebreak

\section{Algebraic Preliminaries}

\noindent We quickly review some algebraic background; see \cite{Eis:alg} or \cite{sch} for additional details.
\begin{defn} Let $I, J$ be ideals in a ring $R$. Then the ideal
    quotient (or colon ideal) is
   \[
     I:J = \{f \in R \mid f\cdot j \in I \mbox{ for all } j \in J\}
   \]
 \end{defn}
 Let $R = \mathbb{C}[z_1,\ldots, z_n]$, and let  ${\mathfrak m} =
 \langle z_1, \ldots z_n\rangle$  denote the homogeneous maximal
 ideal. Then for homogeneous ideals $I^+$ and $\langle g \rangle$
 generated in degree $m$, with $g \not\in I^+$, the condition of \S 1:
  \begin{equation}\label{Contain}
  \langle g \rangle_{m+1} \subseteq I^+_{m+1}
\end{equation}
holds if and only if
\begin{equation}\label{colonInterpretation}
  I^+:g = {\mathfrak m}.
\end{equation}
\noindent Algebraically, this means that ${\mathfrak m}$ is an {\em associated prime} of $I^+$.
\begin{defn}
For an $R$-module $M$, a prime ideal $P$ is an associated prime of $M$
if $P=\ann(m)$ for some $m \in M$. Write $\Ass(M)$ for the
set of associated primes.
\end{defn}
\begin{example}
 Let $R=\mathbb{C}[x,y]$ and $M = R/\langle x^2,xy\rangle$. Then the
 annihilator of $x \in M$ is the ideal $\langle x,y \rangle$, and the
 annihilator of $y \in M$ is the ideal $\langle x \rangle$. Recall that an ideal $Q$ is primary if $xy \in Q$ implies $x \in Q$ or $y^n \in Q$ for some $n \in \Z$. Polynomial rings over a field are Noetherian, and in a Noetherian ring, every ideal $I$ has a primary decomposition
\[
  I = \bigcap\limits_{i=1}^m Q_i, \mbox{ where }Q_i \mbox{ is a primary
    ideal},
\]
and the associated primes of $M$ are defined as the $P_i$ such that $P_i = \sqrt{Q_i}.$
\end{example}

\noindent The relevance to our problem comes from Lemma 1.3.10 of \cite{sch}.
\begin{lemma}\label{ColonAssPrimes}
  The associated primes satisfy $\Ass(R/I) = \{ \sqrt{I:f}  \mid f \in R\}$.
\end{lemma}

\noindent We need one last definition before moving on to proofs of our results.
\begin{defn}
An ideal $I = \langle f_1, \ldots, f_k\rangle \subseteq R$ with $k$
minimal generators is a {\em complete intersection} if each $f_i$ is not a zero divisor on $R/\langle f_1,
\ldots, f_{i-1}\rangle$. Equivalently, the map
\[
R/\langle f_1,\ldots, f_{i-1}\rangle \stackrel{\cdot f_i}{\longrightarrow} R/\langle f_1,\ldots, f_{i-1}\rangle \mbox{ is an inclusion.}
\]
\end{defn}

\noindent From a geometric standpoint, being a complete intersection means that the locus
$V(f_1,\ldots,f_k)$ where the $f_j$ simultaneously vanish has
codimension equal to $k$.

\section{Proof of the theorem}

\begin{prop}\label{nLT3}
If $n \le 3$ and \eqref{Contain} holds, then $P \ge 3.$
\end{prop}
\begin{proof}
First, if $P=1$, then $I^+$ and $\langle g \rangle$ are principal
ideals, and the only way two principal ideals generated in the same
degree $m$ can be equal in degree $m+1$ is if the ideals are themselves
equal, which contradicts our hypothesis. Note that this covers the
case $n=2$, since if $g \not \in I^+$ and \eqref{Contain} holds,
$P\ge 2$.

If $P=2$ (regardless of the number of variables), then either
\begin{enumerate}
  \item{ Case 1: the generators $f_1$ and $f_2$ of $I^+$ share a common facto, say $f_1 =
      h_1F$ and $f_2=h_2F$, with $h_i$ having no common factor.}
  \item{ Case 2: $f_1$ and $f_2$ have no common factor.}
    \end{enumerate}

The ideals $H=\langle h_1,h_2\rangle$ (in Case 1), and $F=\langle
f_1,f_2\rangle$ (in Case 2) are by definition complete
intersections. By Corollary 18.14 of \cite{Eis:alg}, all
associated primes of $H$ and $F$ are of codimension two. Again, this
fact does not depend on the number of variables. Now suppose $n \ge 3$.
Then \eqref{Contain} can only hold in Case 2 if $P \ge 3$. In Case
1, the associated primes of $I^+$ are the union of the associated primes
of $H$ (which are codimension two), along with the codimension one
primes corresponding to irreducible factors of $F$. Thus in Case 1 we
again have that all associated primes are of codimension at most two,
so \eqref{Contain} can only hold if $P \ge 3$.
\end{proof}

In order to show that this result is optimal, we need a bit more
homological algebra. The result we need is a classical
theorem of Bruns \cite{bruns} from 1976.
\begin{theorem}\label{Bruns}
 Given a free resolution $F_\bullet$ with differentials $d_i$, there exists a three-generated
 homogeneous ideal $B$ such that the free resolution $F(B)_\bullet$
 and differentials $d^B_i$  satisfy
 \[
   F_i = F(B)_i \mbox{ and } d_i = d^B_i \mbox{ for all } i \ge 2.
 \]
\end{theorem}
\noindent We will apply Theorem~\ref{Bruns} to the {\em Koszul complex}:
\begin{defn}
  For ${\bf f}=\{f_1, \ldots, f_k\} \subseteq R$, the Koszul complex is
  \[
    K({\bf f}): 0 \longrightarrow R^1 \stackrel{d_k}{\longrightarrow}
    R^k \stackrel{d_{k-1}}{\longrightarrow}  \Lambda^2(R^k)
    \stackrel{d_{k-2}}{\longrightarrow}
    \cdots
  \]
  with differential
  \[
    d_i(e_{j_1} \wedge \cdots \wedge e_{j_i})= \sum\limits_{k=1}^i
    (-1)^kf_{j_k} e_{j_1} \wedge \cdots \widehat{e_{j_k}} \cdots \wedge
    e_{j_i}.
  \]

  \noindent A computation shows that
  \[
    d_i \circ d_{i+1} = 0,
  \]
 so $K({\bf f})$ is a complex; it is exact if and only if $\{f_1, \ldots, f_k\}$ is a complete intersection.
\end{defn}

\begin{prop}\label{nGE4}
If $n \ge 4$, then \eqref{Contain} can hold with $P < n.$
\end{prop}
\begin{proof}
By Lemma~\ref{ColonAssPrimes}, when \eqref{colonInterpretation}
holds, then $g$ satisfies $I^+:g = {\mathfrak m}$ and
${\mathfrak m}$ is an associated prime of $I^+$, since it annihilates
$g$. (Recall that $g$ is a nonzero element of $R/I^+$.) By the Auslander-Buchsbaum Theorem
(Exercise 19.8 of \cite{Eis:alg}), the projective dimension
of $R/I^+$ (i.e., the length of a minimal free resolution) is $n$. Our goal is to
produce an ideal $I^+$ with this behavior--maximal projective
dimension--but with a small number of generators.

To achieve this aim, we apply Theorem~\ref{Bruns} to the Koszul complex of
$\mathfrak{m}$. The result is a three-generated ideal $I^+$ with
${\mathfrak m}$ an associated prime of $I^+$. Lemma~\ref{ColonAssPrimes}
then yields an element $g$ such that
\[
  \sqrt{I^+:g} = {\mathfrak m}.
\]
For $K(\{z_1,\ldots, z_n\})$, the ideal $I^+$ has three homogeneous
generators of degree $n-2$. By Lemma~\ref{ColonAssPrimes} there is a
homogeneous element $g$ satisfying $I^+:g = \mathfrak{m}$, but $g$ could
be of degree greater than $n-2$. A computation using the
software package {\tt Macaulay2} \cite{danmike} shows that for
$n=4$, Theorem~\ref{Bruns} applied to $K({\mathfrak m})$ yields
an ideal $I^+$ with three quadratic generators and that $I^+:{\mathfrak
  m}$ contains a fourth quadric $q \not \in I^+$.
\end{proof}
\begin{example}
There are different choices for the ideal $I^+$.
Perhaps the simplest is
\[
  I^+ = \langle z_4^2, z_2z_3+z_1z_4, z_2^2+z_2z_4\rangle, \mbox{ and a quick check shows } z_2^2 \in I^+: \langle z_1,z_2,z_3,z_4\rangle.
\]
Thus, $g=z_2^2$ gives a minimal example where \eqref{Contain} holds with $P < n$.
\end{example}

\subsection{Directions for future work} A complete answer to Question \ref{ques:min P for signature P,1} will only allow us to prove the Sum of Squares Conjecture when $r(z,\bar{z})$ has signature $(P,N)$ for $N$ small.  To extend this range requires answering the analogue of Question 1 when $N>1$, which is the subject of our ongoing work.






\bibliographystyle{alpha}

\bibliography{references}

\end{document}